\documentclass[conference]{IEEEtran}
\IEEEoverridecommandlockouts
\ifCLASSINFOpdf
\usepackage[pdftex]{graphicx}
\usepackage{subfig}
\else
\fi
\usepackage{amsmath}
\usepackage{amsfonts}

\usepackage[ruled, norelsize]{algorithm2e}
\newtheorem{theorem}{Theorem}

\begin{document}
%
\title{Accelerated Convergent \\Motion Compensated Image Reconstruction}

\author{
\IEEEauthorblockN{Claire Delplancke\thanks{At the time of this work, CD was as the Department of Mathematical Science, University of Bath.}}
\IEEEauthorblockA{EDF Lab Paris Saclay}
\and
\IEEEauthorblockN{Kris Thielemans, \textit{Senior Member, IEEE}}
\IEEEauthorblockA{Institute of Nuclear Medicine\\University College London}
\and
\IEEEauthorblockN{Matthias J.\@ Ehrhardt}
\IEEEauthorblockA{Institute for Mathematical Innovation\\University of Bath}
}


\maketitle

\begin{abstract}
Motion correction aims to prevent motion artefacts which may be caused by respiration, heartbeat, or head movements for example. In a preliminary step, the measured data is divided in gates corresponding to motion states, and displacement maps from a reference state to each motion state are estimated. One common technique to perform motion correction is the motion compensated image reconstruction framework, where the displacement maps are integrated into the forward model corresponding to gated data. For standard algorithms, the computational cost per iteration increases linearly with the number of gates. In order to accelerate the reconstruction, we propose the use of a randomized and convergent algorithm whose per iteration computational cost scales constantly with the number of gates. We show improvement on theoretical rates of convergence and observe the predicted speed-up on two synthetic datasets corresponding to rigid and non-rigid motion. 
\end{abstract}

\section{Introduction}
In standard image reconstruction without Motion Compensation (MC), partial updates on subsets are the cornerstone of acceleration, e.g.\@ OSEM, Block-ART etc, and are used in the provenly convergent Stochastic Primal-Dual Hybrid Gradient algorithm (SPDHG) \cite{PET_SPDHG}, a randomized version of the Primal-Dual Hybrid Gradient algorithm (PDHG) \cite{PDHG}. In the MC framework, an OSEM-like algorithm with partial updates on gates has been presented in \cite{GateOSEM}. Here, we propose to use SPDHG for MC image reconstruction by randomly sampling one gate at each iteration. We derive the theoretical linear convergence rates of PDHG and SPDHG in a strongly convex framework and show that SPDHG's linear convergence improves on PDHG's if the condition number associated to the problem is sufficiently large. We present the results of two experiments on synthetic data corresponding to a single slice of a CT scan of a walnut and of the chest for rigid and non-rigid motion, which highlight the acceleration brought by randomization.

\section{Method}
\subsection{Model and algorithm}
Let us denote by $N$ the number of gates, and model the observed data of gate $i$ by $d_i = A_ix + \epsilon_i$ where $\epsilon_i \sim \mathcal{N}(0,\sigma^2N^{-1})$ and $A_i= A D_i$ is the composition of the displacement operator $D_i: X\rightarrow X$ and the forward operator $A:X \rightarrow Y$, for example the ray transform. To recover the ground-truth image $x$, we solve the variational model:
\begin{align*}
\underset{x}{\min}\quad O(x) &:= g(x)+\sum_{i=1}^Nf_i(A_ix),
\end{align*}
where $f_i(x) = N^{-1} \|A_ix - d_i\|^2$, $g(x)=\alpha \|x\|^2$ and the division by $N$ is used for appropriate scaling of the variance with increased number of gates. By convexity and regularity of the $f_i$'s and $g$, this minimization problem can be rewritten as a saddle point problem and solved by Algorithm \ref{alg:AlgoScheme}. In PDHG, all dual variables $(y_i)$ are updated at each iteration, corresponding to the choice $p_i=1$ for all $1\leq i \leq N$, whereas in SPDHG only one dual variable is randomly selected with probability $p_i=1/N$ and updated. An \textit{epoch} corresponds to one iteration of PDHG and $N$ iterations of SPDHG, amounting to the same number of forward and adjoint operations.
\vspace{-0.01\textheight}
\begin{algorithm}
  \SetKwInOut{Input}{Input}
  \SetAlgoLined
  \Input{$(\sigma_i),\,\tau,\,\theta,\,(p_i)$. \textbf{Init.\@:} $x=z=\bar{z}=0$, $y=0$.}
  \For{$k=0, \dots,$}{
      (i) $x^{k+1} = \text{prox}_{\tau g}(x^k - \tau z^k)$\;
      (ii) Draw $\mathbb{S} \subset\left\{ 1, \dots N\right\}$ s.\@t.\@ $\mathbb{P}(i\in \mathbb{S})=p_i$\;
      $y_i^{k+1} = \begin{cases} \text{prox}_{\sigma_i f_i^*}(y_i^k + \sigma_i A_i x^{k+1}) & \text{if } i\in \mathbb{S}\\
      							y_i^k & \text{else}.\end{cases}$\\
      (iii) $z^{k+1} = z^k + \sum_{i\in \mathbb{S}}A_i^*(y_i^{k+1}-y_i^k)$ \;
      $\bar{z}^{k+1} = \bar{z}^k + \sum_{i\in \mathbb{S}}(1+\theta/p_i) A_i^*(y_i^{k+1}-y_i^k)$.
  }
  \caption{SPDHG} \label{alg:AlgoScheme}
\end{algorithm}
\vspace{-0.01\textheight}
\subsection{Rates of convergence}
For strongly convex $(f_i^*)$ and $g$, as is the case in our model, and admissible step-sizes $(\sigma_i),\tau,\theta$, the PDHG and SPDHG converge linearly to the saddle point $(x^*,y^*)$ \cite{SPDHG}: at epoch $K$,
\begin{align*}
\|(x_{\text{PDHG}}^K,y_{\text{PDHG}}^K) - (x^*,y^*)\|^2 &\leq C (r^{\text{PDHG}}_N)^K\\
\mathbb{E}\left[\|(x_{\text{SPDHG}}^K,y_{\text{SPDHG}}^K) - (x^*,y^*)\|^2\right] &\leq \tilde{C} (r^{\text{SPDHG}}_N)^K.
\end{align*}
Let $l(\kappa,n)=\left(1-2n^{-1}(1+\sqrt{1+\kappa})^{-1} \right)^n$. By applying \cite[Example 3]{SPDHG} to the current optimization problem, we find that:
\begin{theorem}
For $N$ gates and  well-chosen step-sizes $(\sigma_i),\tau,\theta$, the condition number and per-epoch convergence rate of SPDHG and PDHG are
\begin{align*}
\kappa^{\text{SPDHG}}_N &= \max_i \|A_i\|^2/(\alpha N), & \kappa^{\text{PDHG}}_N &= \|(A_0,\dots,A_N)\|^2/(\alpha N),\\
r^{\text{SPDHG}}_N &=l(\kappa^{\text{SPDHG}}_N,N), &r^{\text{PDHG}}_N &=l(\kappa^{\text{PDHG}}_N,1).
\end{align*}
\end{theorem}
To further exploit these formulas, let us denote $\kappa=\|A\|^2/\alpha$ and notice that the displacement map $D_i$ has usually little influence on the norm of $A_i=AD_i$, so that $\max_{1\leq i\leq N}\|A_i\|^2\approx \|A\|^2$ and $\|(A_0,\dots,A_N)\|^2\approx N\|A\|^2$. For example, these approximations are satisfied with $0.03$ and $0.08$ precision respectively in the numerical applications below. It follows that $r^{\text{SPDHG}}_N \approx l(\kappa N^{-1},N)$ and $r^{\text{PDHG}}_N \approx l(\kappa,1)$. We show that for $\kappa\geq 10$ (a moderately conditioned problem) $ l(\kappa N^{-1},N) < l(\kappa,1)$, hence SPDHG converges at a faster rate than PDHG.

\begin{figure}[htb]
\centering
\captionsetup[subfigure]{indention=.025\textwidth}
\subfloat[Rigid motion]{\includegraphics[width=.22\textwidth]{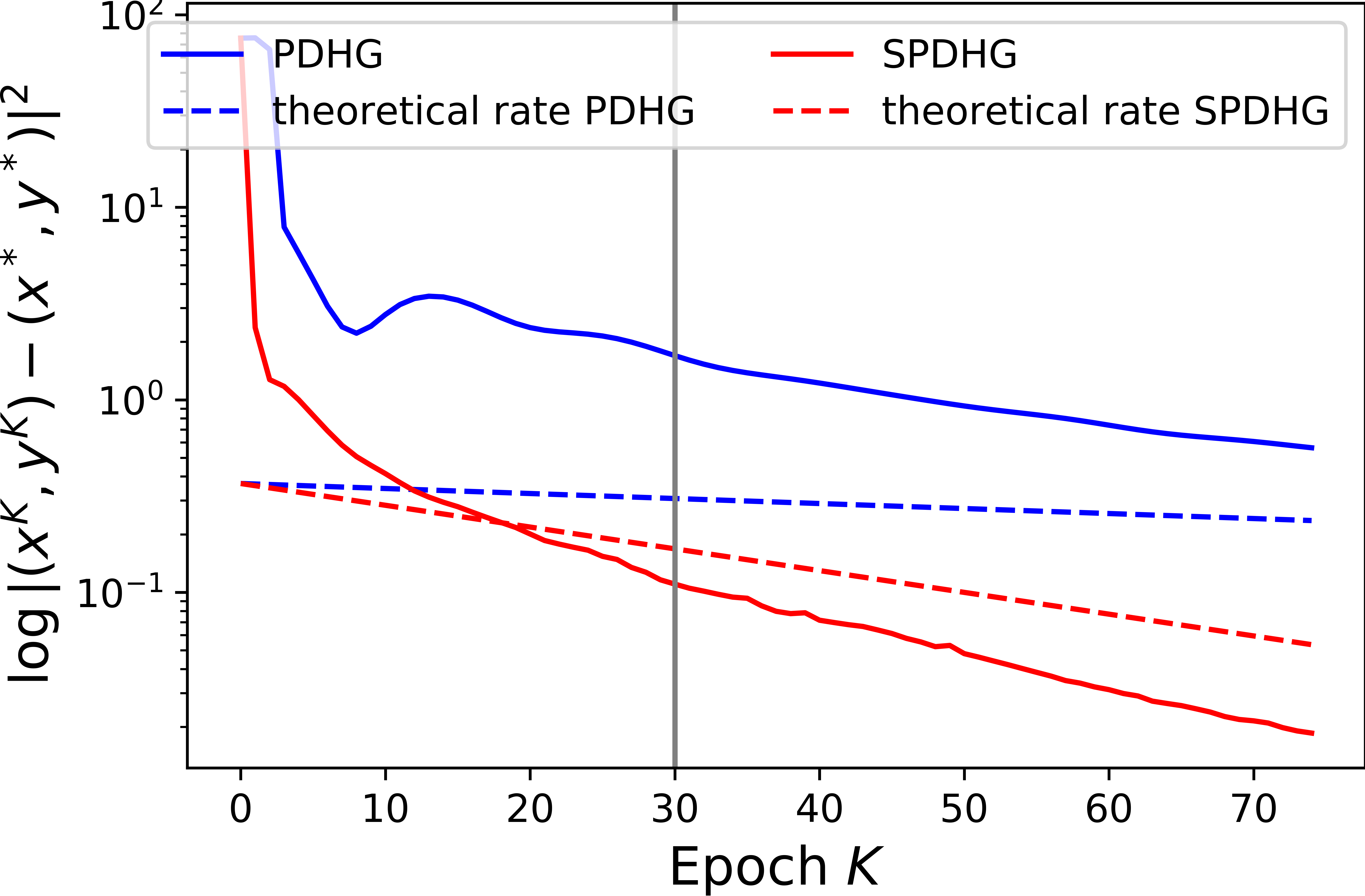}}\qquad
\subfloat[Non-rigid motion]{\includegraphics[width=.22\textwidth]{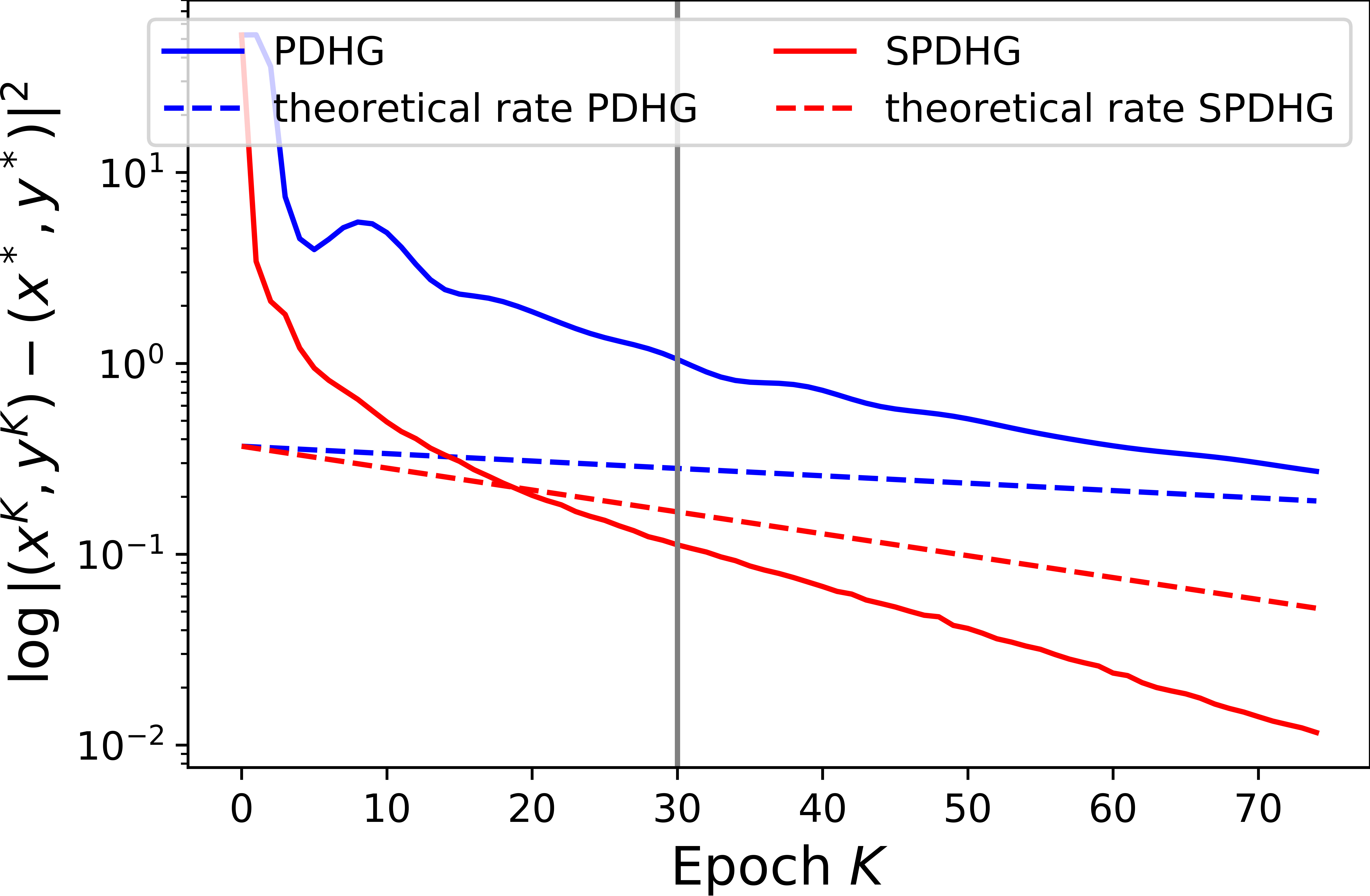}}
\caption{SPDHG's linear convergence is faster than PDHG's}
\label{fig:meas}
\end{figure}

\begin{figure}[htb]
\centering
\captionsetup[subfigure]{indention=.025\textwidth}
\subfloat[First motion state]{\includegraphics[width=.15\textwidth]{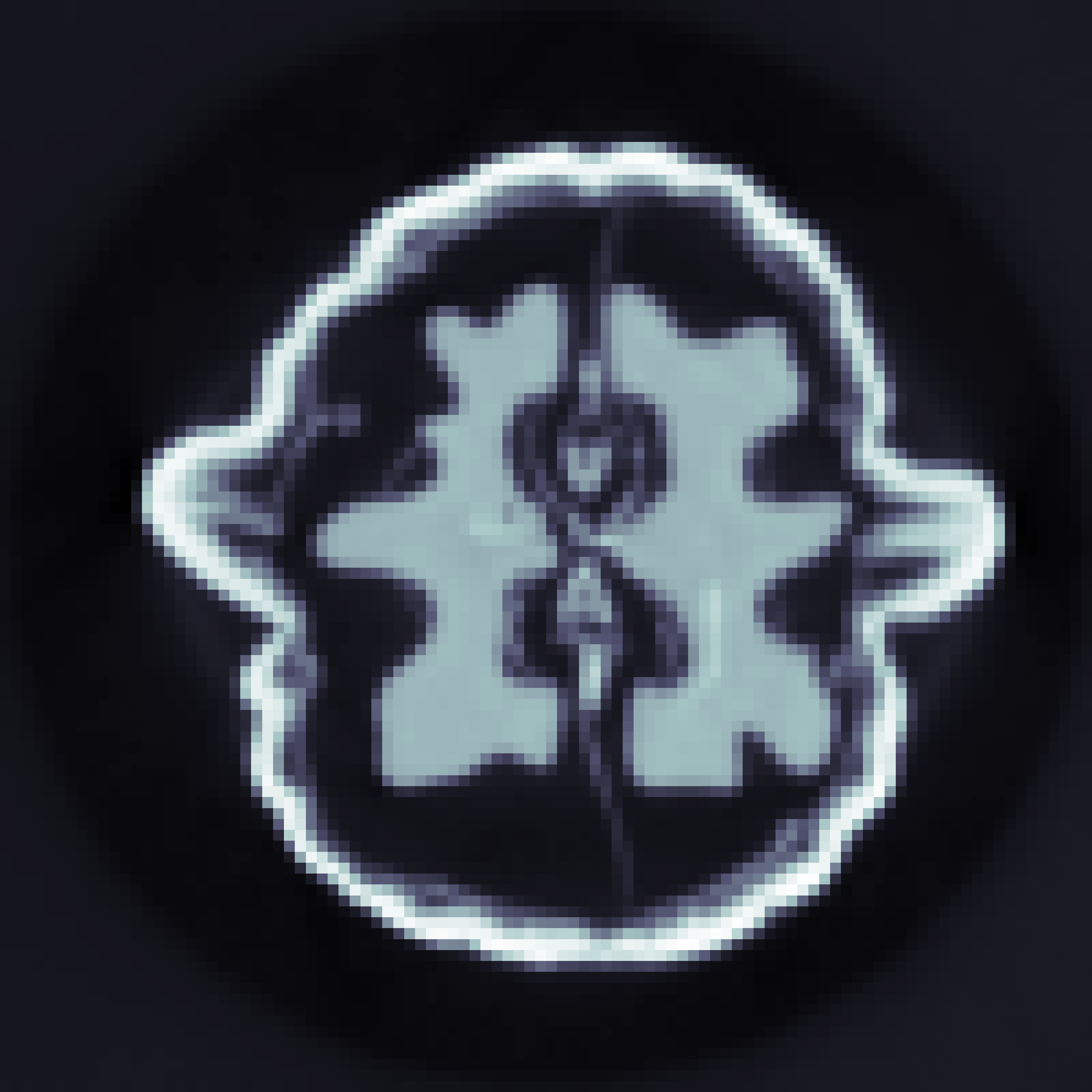}}\,
\subfloat[Last motion state]{\includegraphics[width=.15\textwidth]{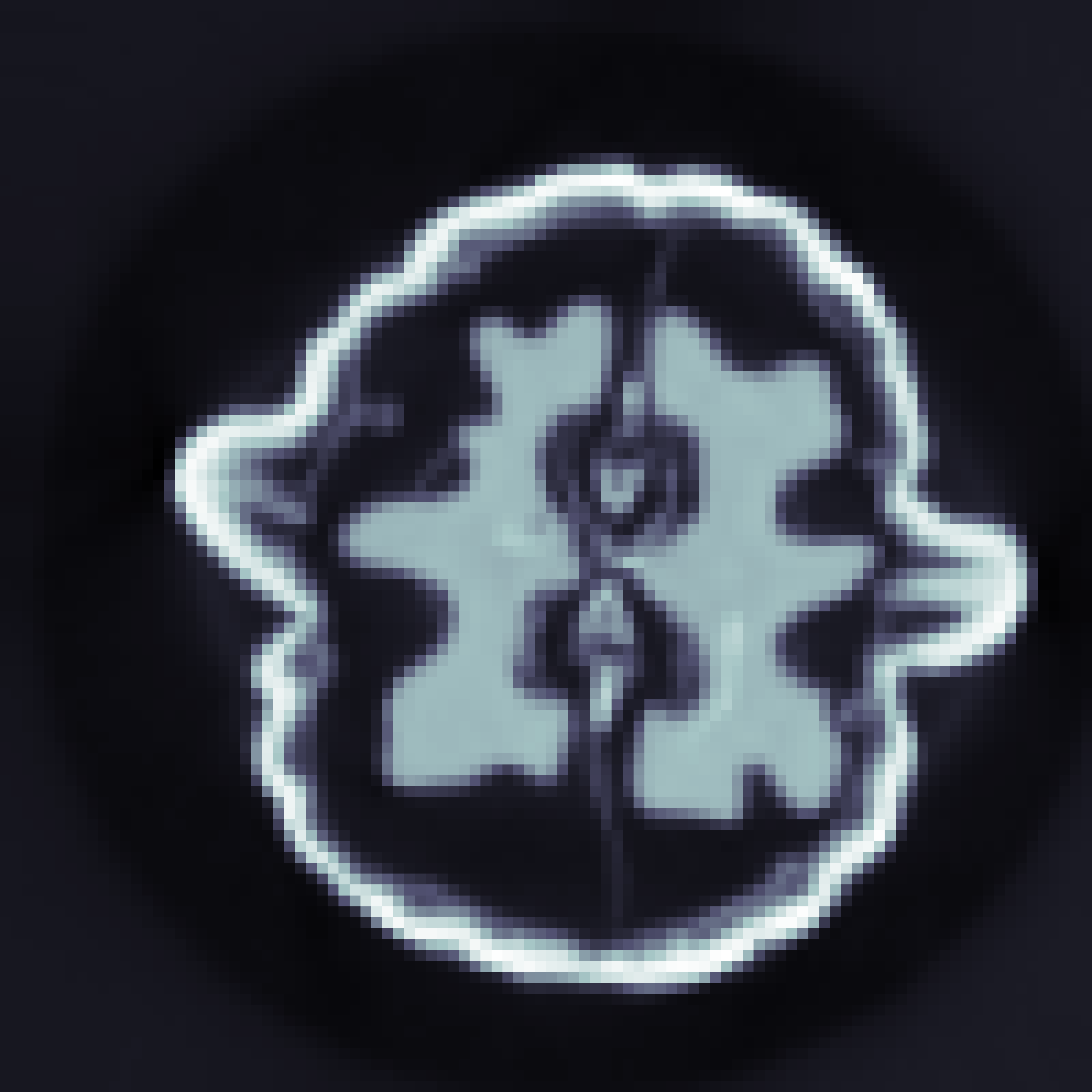}}\,
\subfloat[Converged no-MC]{\includegraphics[width=.15\textwidth]{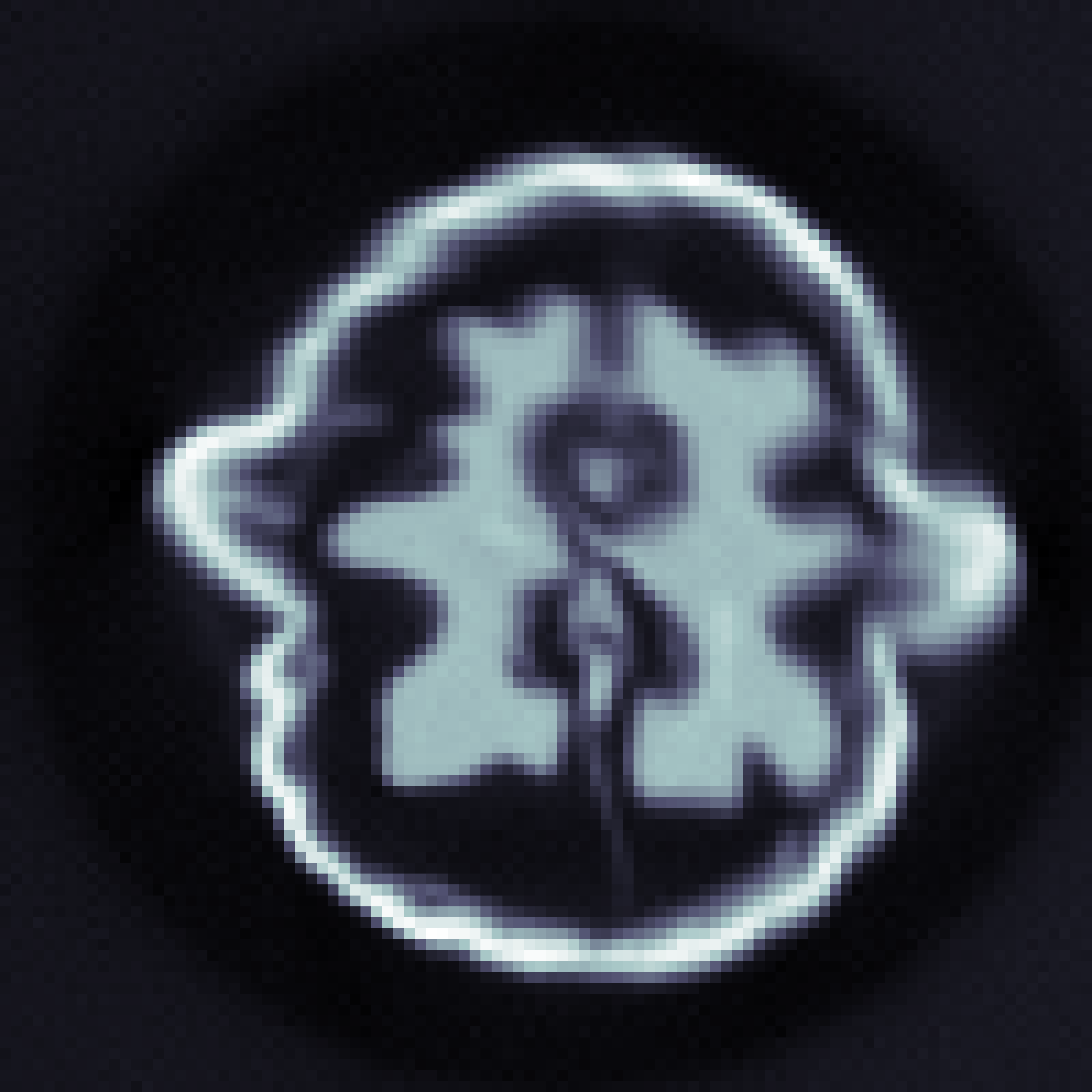}}\\
\subfloat[Converged MC]{\includegraphics[width=.15\textwidth]{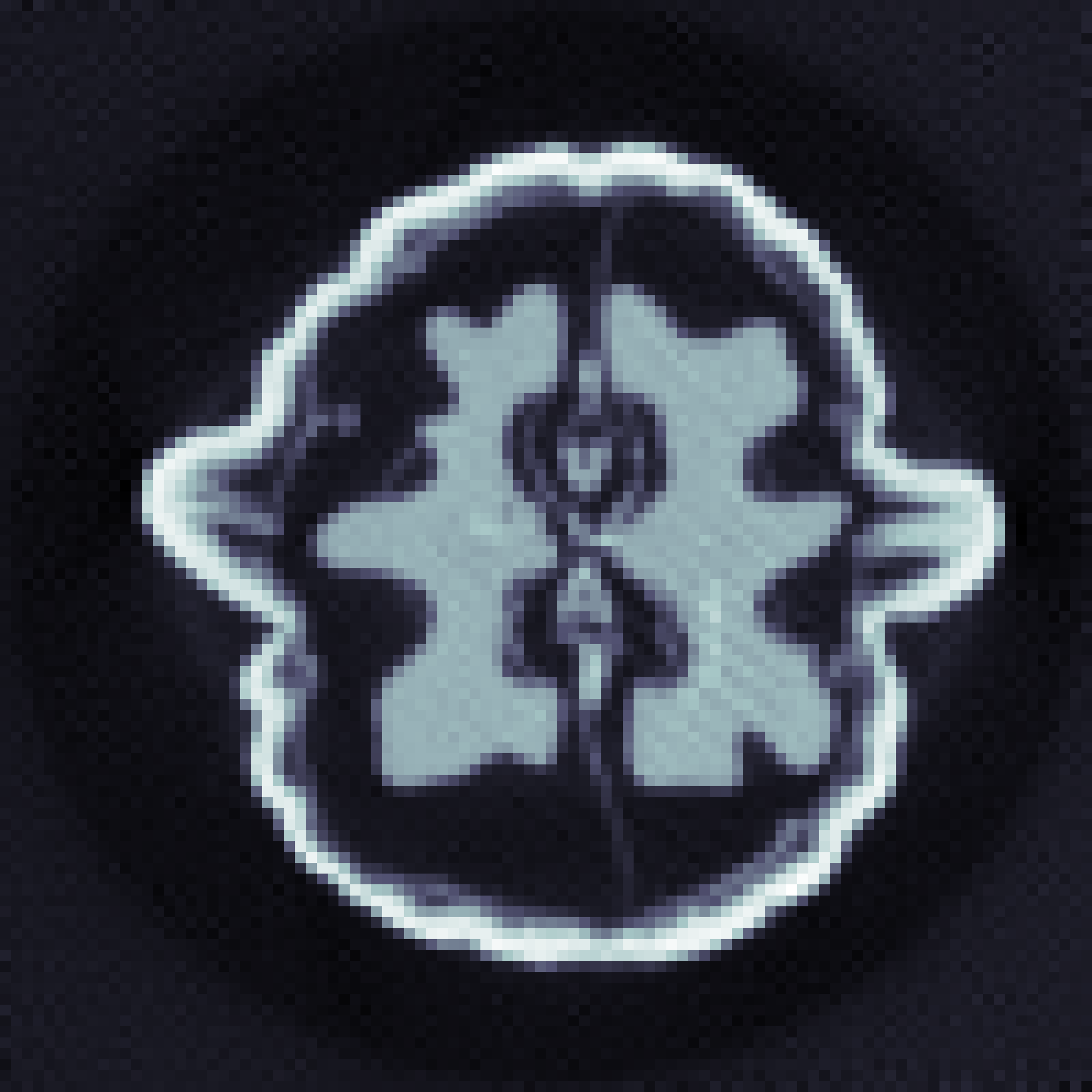}}\,
\subfloat[MC SPDHG after 30 epochs]{\includegraphics[width=.15\textwidth]{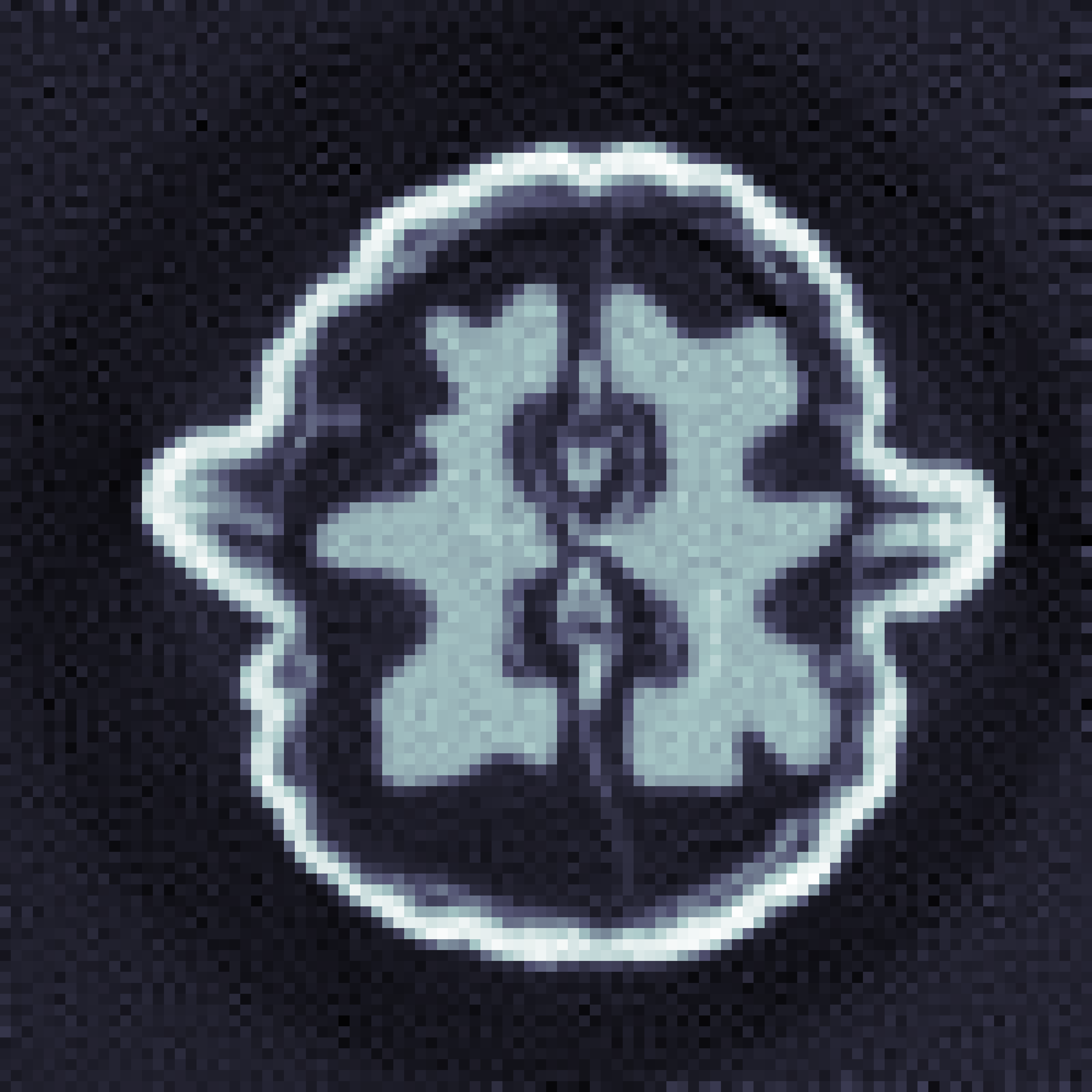}}\,
\subfloat[MC PDHG after 30 epochs]{\includegraphics[width=.15\textwidth]{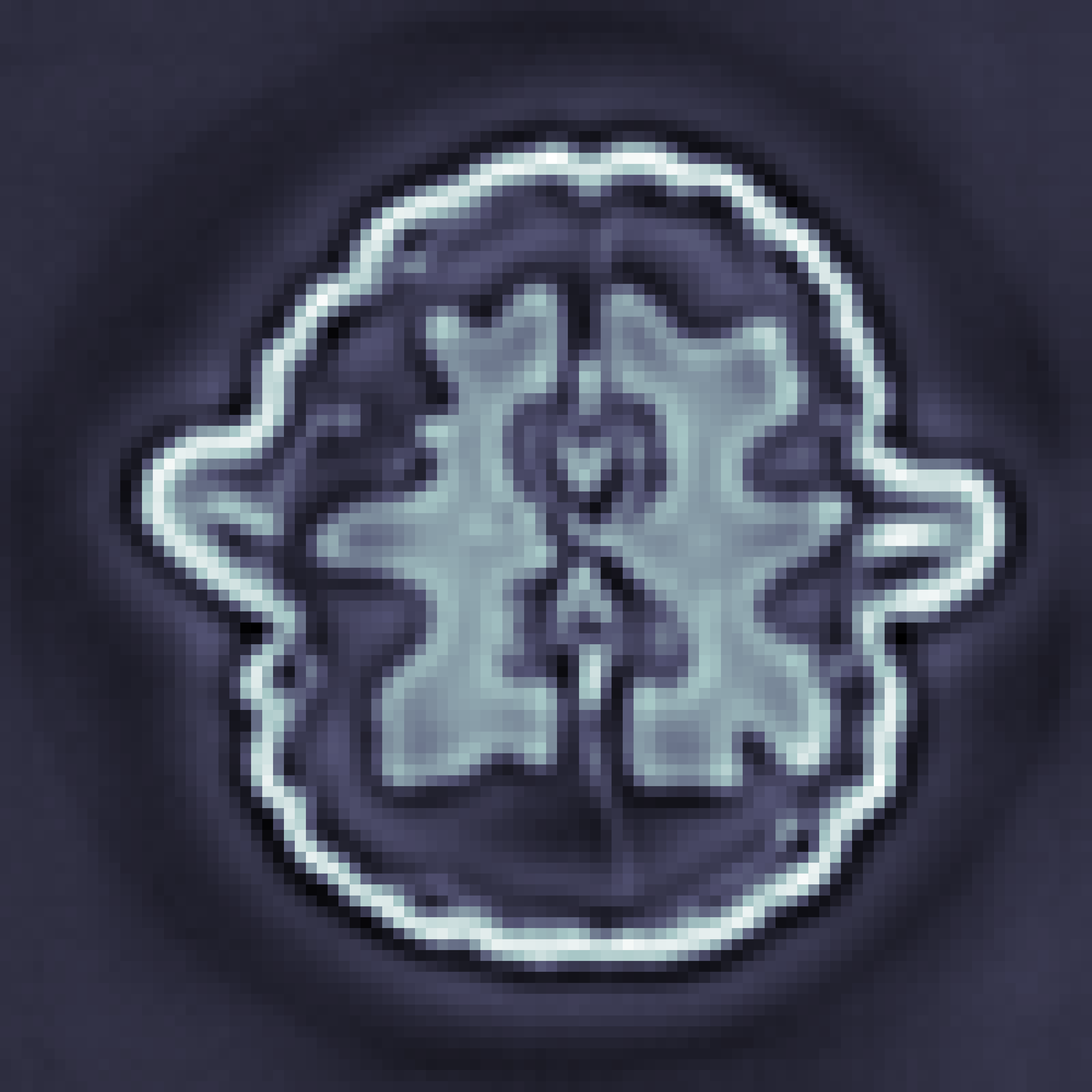}}
\caption{Rigid motion}\label{fig:rigid}
\end{figure}

\begin{figure}[htb]
\centering
\captionsetup[subfigure]{indention=.025\textwidth}
\subfloat[First motion state]{\includegraphics[width=.15\textwidth]{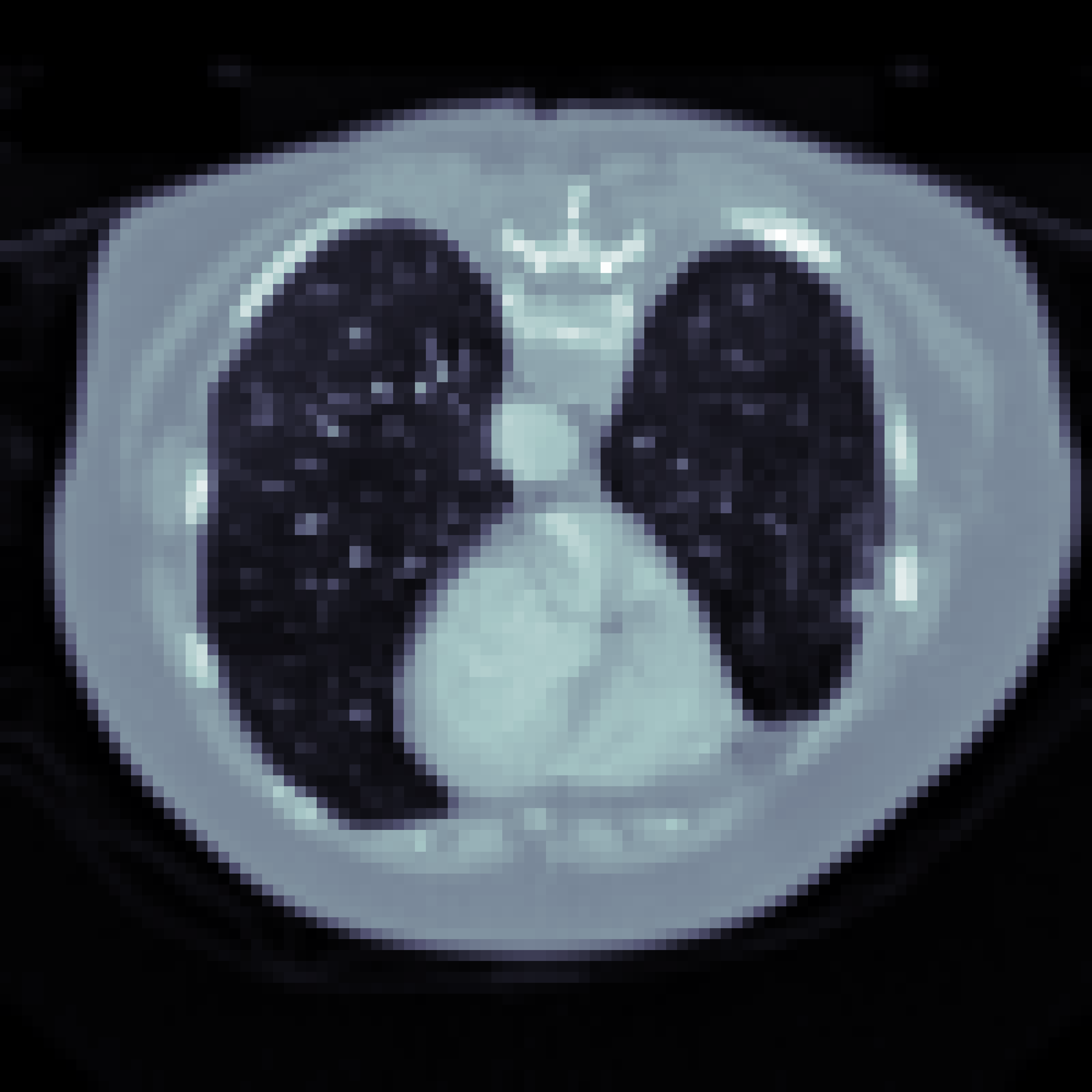}}\,
\subfloat[Last motion state]{\includegraphics[width=.15\textwidth]{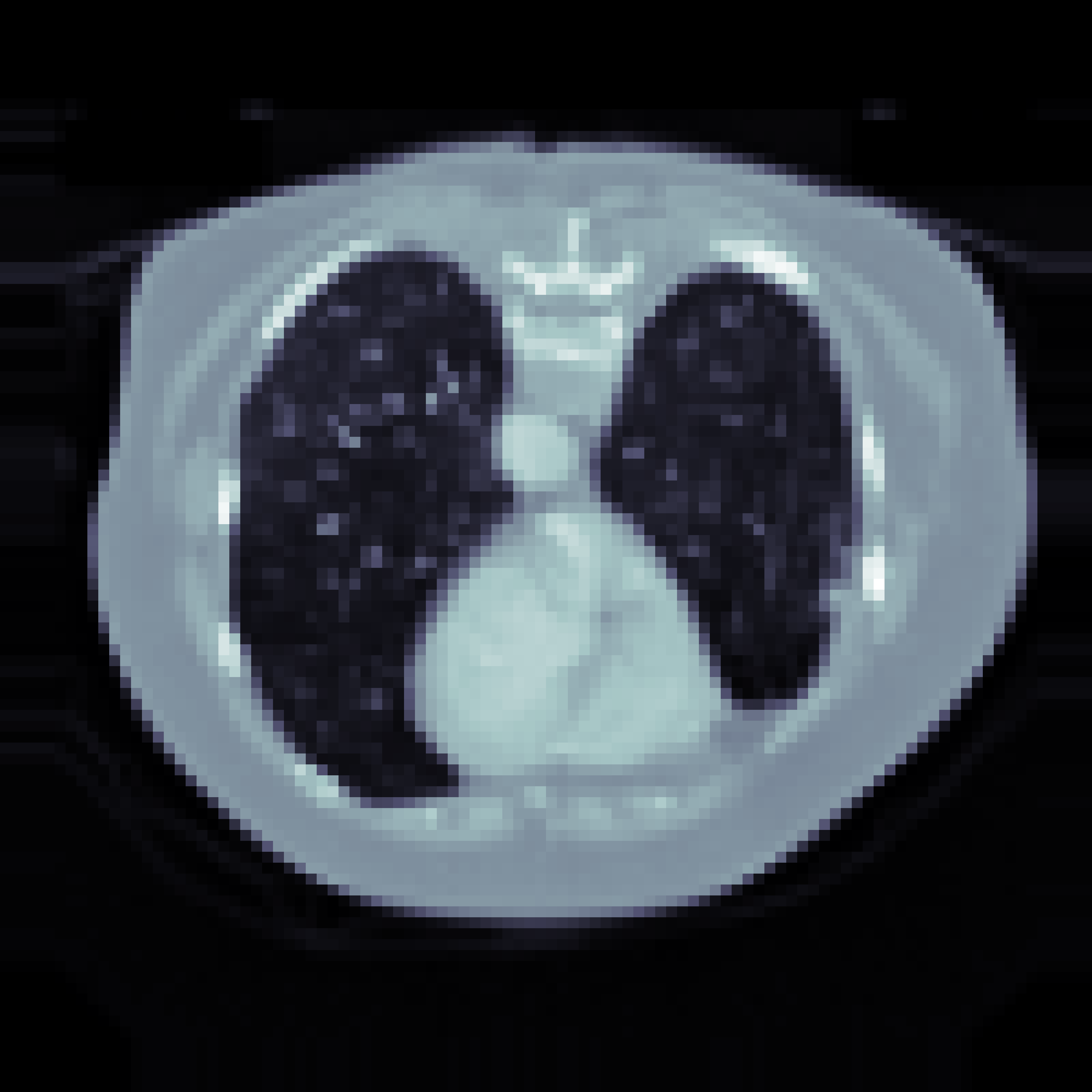}}\,
\subfloat[Converged no-MC]{\includegraphics[width=.15\textwidth]{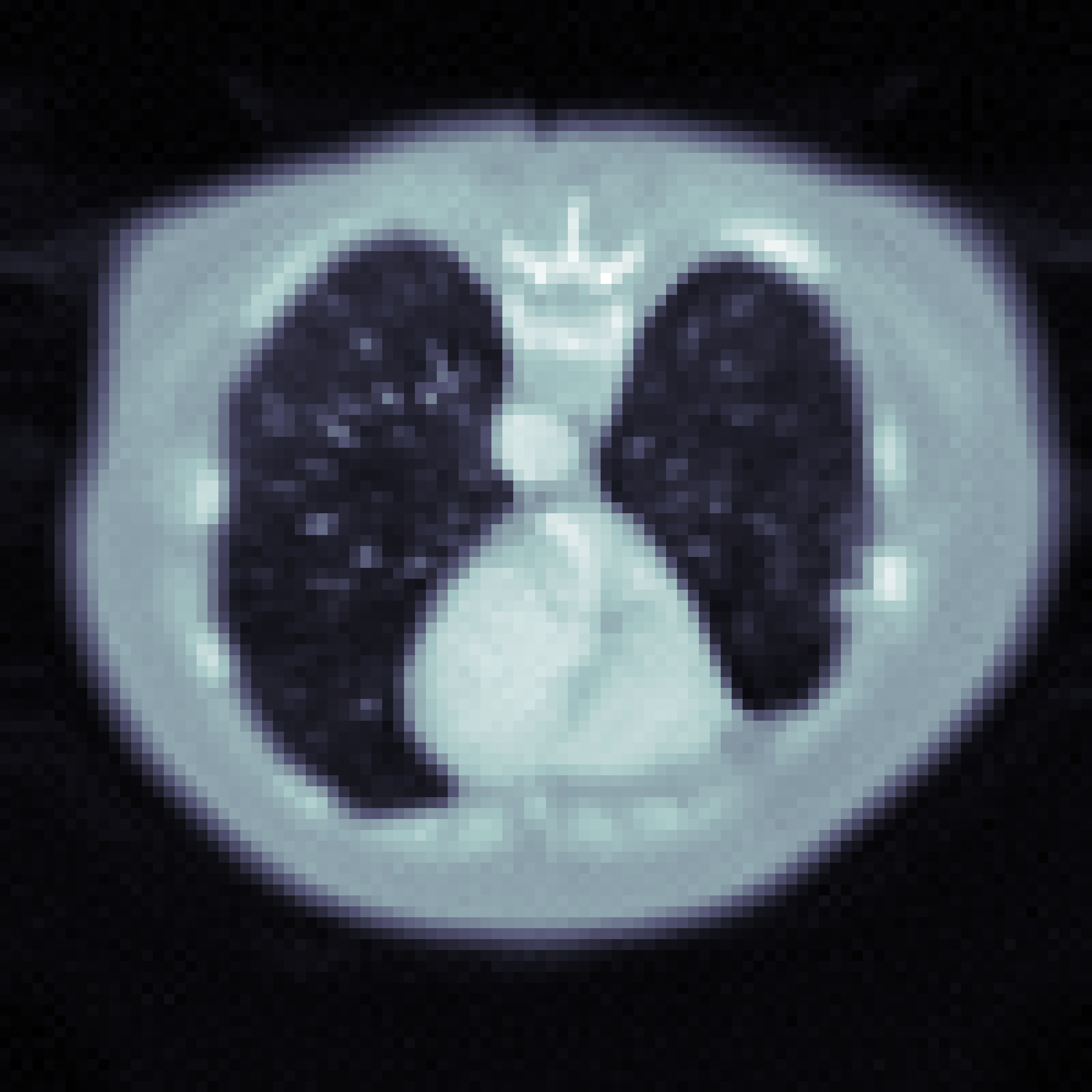}}\\
\subfloat[Converged MC]{\includegraphics[width=.15\textwidth]{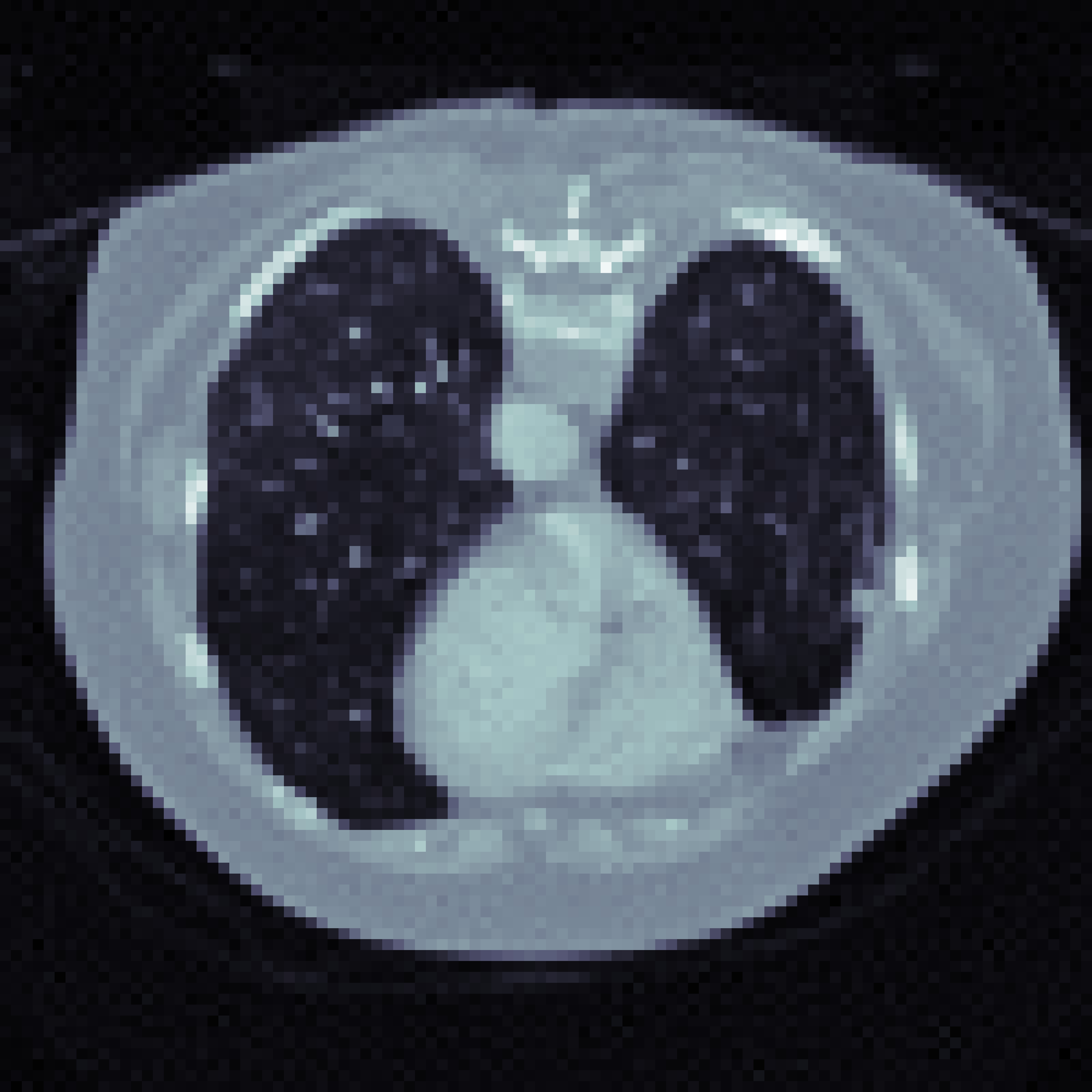}}\,
\subfloat[MC SPDHG after 30 epochs]{\includegraphics[width=.15\textwidth]{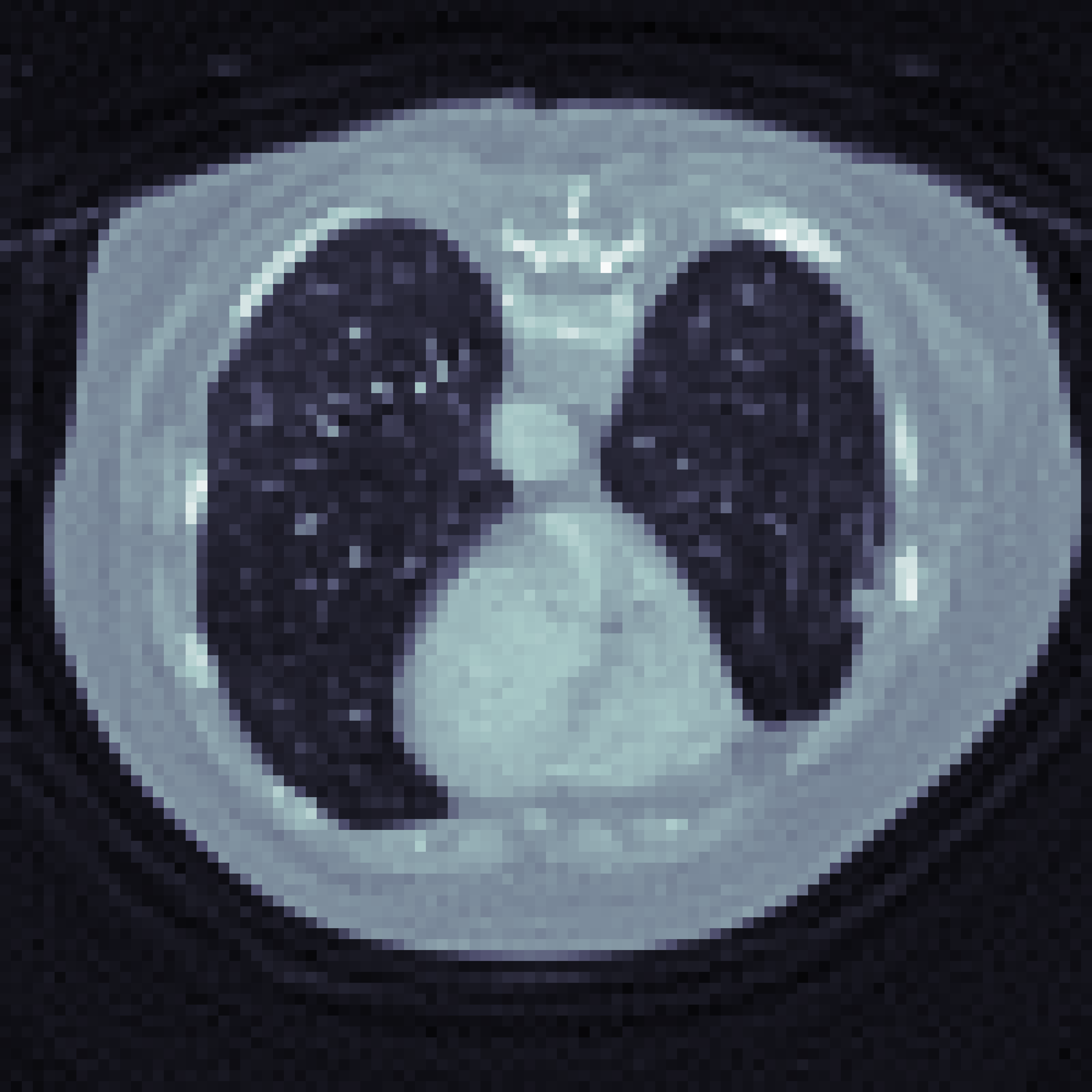}}\,
\subfloat[MC PDHG after 30 epochs]{\includegraphics[width=.15\textwidth]{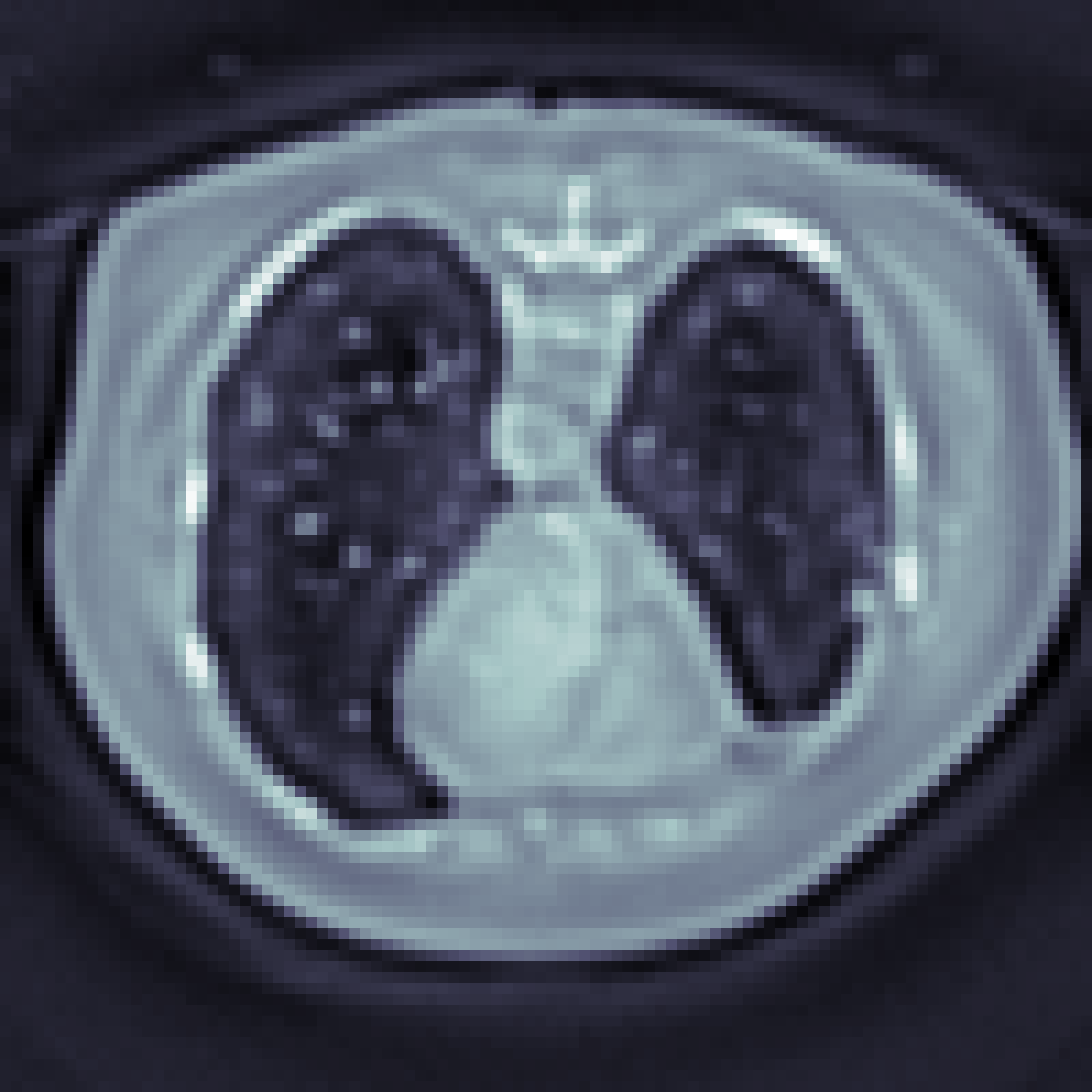}}
\caption{Non-rigid motion}\label{fig:nonrigid}
\end{figure}

\section{Numerical experiments}
\textit{Description.} 
The numerical experiments use the open-source package Astra-Toolbox \cite{AstraToolbox} to compute the ray transform $A:\mathbb{R}^{100\times100}\rightarrow \mathbb{R}^{200\times200}$ in a 2D parallel geometry with $200$ angles and the open-source package ODL \cite{ODL} for all other operations, including the proximal operators. We consider two synthetic datasets. The first one corresponds to a 2D CT image of a walnut \cite{CTwalnut} subject to rigid motion (rotation and translation) for $N=20$ gates, as shown in Figure \ref{fig:rigid}, subplots (a) and (b). The second one corresponds to a 2D slice of a chest CT scan \cite{CTchest} subject to non-rigid motion (linear dilatation) for $N=10$ gates, as shown in Figure \ref{fig:nonrigid}, subplots (a) and (b). For both models, the condition number defined above is $\kappa \approx 70$.

\textit{Results.}
As expected, the converged no-MC reconstructions shown in Figures \ref{fig:rigid} and \ref{fig:nonrigid}, subplot (c), look blurry in comparison to the converged MC reconstructions shown in Figures \ref{fig:rigid} and \ref{fig:nonrigid}, subplot (d).

We observe in Figure \ref{fig:meas} that SPDHG and PDHG converge linearly with respect to the $l_2^2$-distance. We picked arbitrary values $C,\,\tilde{C}$ for the theoretical rates in order to make the comparison between the slopes visually easy. The observed convergence rates are consistent with the theoretical ones, with SPDHG faster than PDHG.

The speed gain is confirmed when comparing SPDHG and PDHG reconstructions both at epoch $30$ for rigid (Figure \ref{fig:rigid}) and non-rigid (Figure \ref{fig:nonrigid}) motion: the SPDHG reconstruction in subplot (e) appears visually closer to the optimal MC reconstruction in subplot (d) than the PDHG reconstruction in subplot (f).

\section{Conclusions}
We presented how to accelerate MC image reconstruction by randomly sampling over the gates with SPDHG. In order to further improve the speed-up gain, our next step is to investigate the performance of SPDHG when sampling over both gates and over data subsets (e.g. angles).

\bibliographystyle{IEEEtran}
\bibliography{IEEEabrv,biblio_abb.bib}
\end{document}